\documentclass[12pt,twoside]{article}
\usepackage{amssymb}
\usepackage{amsthm}
\usepackage{graphicx}
\def\R{{\mathbb R}}
\def\N{{\mathbb N}}
\def\C{{\mathbb C}}
\def\Q{{\mathbb Q}}
\def\T{{\mathbb T}}
\def\Z{{\mathbb Z}}
\def\kasten{$~~\mbox{\hfil\vrule height6pt width5pt depth-1pt}$ }
\newtheorem{Theorem}{Theorem}[section]

\newtheorem{Definition}[Theorem]{Definition}

\newtheorem{Remark}[Theorem]{Remark}

\begin{document}
\par\noindent
{\Large \bf Partly Divisible Probability Measures
\newline
on Locally Compact Abelian Groups}

\bigskip
\par\noindent
by Sergio Albeverio of Bonn, Hanno Gottschalk of Roma and

Jiang-Lun Wu of Bochum

\bigskip\par
\begin{abstract}{A notion of admissible probability measures
$\mu$ on a locally compact Abelian group (LCA-group) $G$ with
connected dual group $\hat G=\R^d\times \T^n$ is defined. To
such a measure $\mu$, a closed semigroup
$\Lambda(\mu)\subseteq (0,\infty)$ can be associated, such that,
for $t\in \Lambda(\mu)$, the Fourier transform to the power $t$,
$(\hat \mu)^t$, is a characteristic function. We prove that the
existence of roots for non admissible probability measures
underlies some restrictions, which do not hold in the admissible
case. As we show for the example $\Z_2$, in the case of
LCA-groups with non connected dual group, there is no canonical
definition of the set $\Lambda(\mu)$.}
\end{abstract}
\bigskip

\section{Introduction}

Let $(G,+)$ be a locally compact Abelian group (LCA-group) and
$X$ a $G$-valued random variable (measurable ~with respect to
the Borel $\sigma$-algebra on $G$) on an arbitrarily given
probability space with probability measure $P$. As in the real
random variable case, it is interesting to know whether $X$ can
be decomposed into a sum of $n$ independent identically
distributed $G$-valued random variables $X_1,\ldots,X_n$, say,
$$X=X_1+\cdots+X_n , ~~~ P-\mbox{a.e.}.$$
Alternatively, in terms of $\mu = P\circ X^{-1}$, the image
(probability) measure of $P$ under $X$ , the above equation
becomes the following $n$-fold convolution equation
\begin{equation}
\label{1.1eqa}
\mu={\underbrace{\mu_{1\over n}\ast\cdots\ast\mu_{1\over n}}_{n
~ \mbox{\tiny times}} ,}
\end{equation}
where $\mu_{1\over n}=P\circ X_1^{-1}$ is the image (probability)
measure of $P$ under $X_1$. Furthermore, Equation (\ref{1.1eqa})
can be written in terms of characteristic functions (via the
Fourier transform) as follows
$${\hat\mu}=(\hat{\mu_{1\over n}})^n \, . $$
We may generalize this problem and ask whether for $q={m\over n}$
with $n,m\in \N$ being mutually prime (abbr. m.p.) there exists a
characteristic function (c.f.) $\hat{\mu_{m\over n}}$ such that
\begin{equation}
\label{1.2eqa}
(\hat \mu)^m=(\hat{\mu_{m\over n}})^n ~.
\end{equation}
Let ${\cal M}_1(G)$ denote the space of probability measures on
$G$. Given a measure $\mu\in{\cal M}_1(G)$, we may thus define a
set of positive rational numbers $\Lambda^{\mbox{\rm alg}}(\mu)$
as ~follows:
$$\Lambda^{\mbox{\rm alg}}(\mu)=\{ q={m\over n}: n,m\in \N
\mbox{ m.p. and } \exists \mbox{ c.f. } \hat{\mu_{m\over n}}
\mbox{ such that } (\hat\mu)^m=(\hat{\mu_{m\over n}})^n \} .$$
This set is clearly designed to contain the information about
the (algebraic) divisibility of $\mu$.

In the case where $\hat\mu$ is a real function with $0<\hat\mu
\le1$, we can present a candidate for the solutions of Equation
(\ref{1.2eqa}) as follows
$$\hat{\mu_{m\over n}}=(\hat\mu)^{m\over n}$$
which leaves open the question whether ${m\over n}$ is a number
such that $({\hat\mu})^{m\over n}$ is a characteristic function.

So the (general) problem which arises in connection with the
above observation can be roughly formulated as follows: for
which kind of probability measures $\mu$ is the set
$$\Lambda(\mu)=\{t\in(0,\infty): ({\hat\mu})^t ~~ \mbox{is a
characteristic function}\}$$
well defined and what kind of subsets of $(0,\infty)$ do possess
the property of ~being $\Lambda(\mu)$-sets for some probability
measure $\mu$?

We have studied this question, which in the case $G=\R$ goes back
to D. Dugu\'e \cite{Du}, in a previous note for this special case
\cite{AGW}, where also concrete examples have been discussed. More
examples can be found in \cite{BS,Du}. All these examples together
indicate that there exist a great variety of sets $S\subseteq (0,
\infty)$ such that a measure $\mu$ with $\Lambda(\mu)=S$ can be
found.

In this paper we extend these considerations to the framework of
locally compact Abelian groups. In particular, in Section 2 we
give a proper definition of the set $\Lambda(\mu)$ for the case
of a LCA-group with connected dual group $\hat G$. This leads to
the notion of admissible probability measures on $G$. In
Section 3 we show that the algebraic divisibility given by
$\Lambda^{\mbox{alg}}(\mu)$ for non admissible probability
measures $\mu$ is restricted by some rather general
considerations, which do not apply to the admissible case. In
Section 4 we discuss the situation $G=\Z_2$ as the simplest
example of a LCA-group with non connected dual group and we show
that there is no canonical definition of the set $\Lambda(\mu)$
in this case.

\section{Admissible probability measures on LCA-
\newline
groups with connected dual groups}

In this section we extend the notion of admissible probability
distributions on $\R$ given in \cite{AGW} to the case of
probability measures on a LCA-group $(G,+)$ with arcwise
connected dual group $(\hat G,\cdot)$. Recall that ${\cal M}_1(G)$
is the set of probability measures on $G$. We say
$\mu\in{\cal M}_1(G)$ \underline{admissible}, if there is some
continuous function $\psi:\hat G \to \C$ with $\psi (1)=0$ such
that $\hat\mu = e^\psi$. We shall denote the set of all admissible
probability measures on $G$ by ${\cal M}_1^a(G)$.

The function $\psi$ is uniquely determined. In fact, if $\psi_1$
and $\psi_2$ fulfil the above conditions, then
$$\psi_1(\gamma)-\psi_2(\gamma)=2\pi i\, k(\gamma) ~~~\mbox{with}
~ k(\gamma) \in \Z ~~~\forall \gamma \in \hat G .$$
Since the left hand side of the above equation is continuous, the
right hand side must also be continuous and is thus constant,
since a topological space $X$ is connected if and only if every
continuous map from $X$ to $\Z$ is constant. Evaluation of the
above equation at $\gamma=1$ now yields that $\psi_1=\psi_2$.
For $\mu \in {\cal M}_1^a(G)$ we call the function
$\psi:\hat G\to\C$, which is uniquely determined by the conditions
specified above, the \underline{second characteristic} associated
with $\mu$.

We now give the main definition of this article in analogy to
\cite{AGW}:

\begin{Definition}
\label{2.1def}
Let $G$ be a LCA-group with arcwise connected dual group $\hat G$.
For $\mu\in{\cal M}_1^a(G)$ with second characteristic $\psi$ we
define the set $\Lambda(\mu)$ of positive real numbers as
$$\Lambda (\mu) = \{ t>0: e^{t\psi} ~~ \mbox{is a characteristic
function}~\}~.$$
\end{Definition}

\begin{Remark}
\label{2.1rem}
{\rm (i) By Theor\`eme 4 of \cite{Di}, $\hat G$ is arcwise
connected $\Rightarrow\hat G=\R^d\times\T^n\Leftrightarrow G
=\R^d\times \Z^n$, where $\T^n$ stands for the $n$-dimensional
torus. In particular, arcwise connected LCA-groups are also
locally arcwise connected.

\noindent (ii) Let ${\cal L}(G)$ denote the complex measures on
$G$ which have a logarithm in the Banach algebra of complex
measures with the convolution as multiplication and
${\cal L}_1(G)$ the probability measures in ${\cal L}(G)$.
Then ${\cal L}_1(G)\subseteq {\cal M}_1^a(G)$ where equality
holds if $\hat G$ is compact (i.e. $G=\Z^n$) and the inclusion
is proper if $\hat G$ is not compact ($G=\R^d\times\Z^n,
d\not = 0$), since e.g. any Gaussian measure supported on one
of the copies of $\R$ has unbounded second characteristic $\psi$
and thus does not belong to ${\cal L}_1(G)$. For results on the
characterization of ${\cal L}(G)$ see e.g. \cite{GW,Ha,Hi}.

\noindent (iii) $\Lambda(\mu)\subseteq (0,\infty)$ is a closed
semigroup under addition with $\N\subseteq\Lambda(\mu)$. Either
$\Lambda(\mu)=(0,\infty)$ or $\exists \lambda >0$ such that
$\Lambda(\mu)\subseteq[\lambda,\infty)$. If the interior of
$\Lambda(\mu)$ is not empty, then $\exists\lambda >0$ such that
$\Lambda(\mu) \supset [\lambda,\infty)$, cf. Prop. 2.1 and 2.2
of \cite{AGW}.

\noindent (iv) $\mu \in M_1^a(G)$ is infinitely divisible if
and only if $\Lambda(\mu)=(0,\infty)$ (Schoenberg's theorem
\cite{BeFo}). If $\Lambda(\mu)\not = (0,\infty)$, we say $\mu$
is \underline{partly divisible} and if $\Lambda(\mu)=\N$ we call
$\mu$ \underline{minimally divisible} (see \cite{AGW} for
examples). }
\end{Remark}

We now want to give a topological characterization of the Fourier
transform of admissible probability measures.

The assumptions on $\hat G$ allow us to apply the theory of
covering spaces to characterize the Fourier transforms of
admissible probability measures.

Let us first recall some basic notions of algebraic topology
following \cite{Ma}. Let $X,Y$ be topological spaces. Then
$\phi : (X,x_0)\to (Y,y_0)$ denotes a continuous map from $X$
to $Y$ which maps $x_0$ to $y_0$. The fundamental group of $X$
(resp. $Y$) based at $x_0$ (resp. $y_0$) is denoted by
$\pi(X,x_0)$ (resp. $\pi(Y,y_0)$). Then $\phi_*:\pi(X,x_0)\to
\pi(Y,y_0)$ is the group homomorphism induced by $\phi $.

A covering space $(\tilde X, \tilde x_0,\rho)$ of a topological
space $(X,x_0)$ consists of a topological space $(\tilde X,
\tilde x_0)$ and a map $\rho :(\tilde X,\tilde x_0)\to (X,x_0)$
such that for every $x\in X$ there exists an ~arcwise connected
open $x$-~neighborhood $U\subseteq X$ such that each ~arcwise
connected component of $\rho^{-1}(U)$ is homeomorphic to $U$
under the restriction of $\rho$ to this component. A lifting of
the map $\phi:(Y,y_0)\to(X,x_0)$ to $(\tilde X,\tilde x_0,\rho)$
is a map $\tilde\phi :(Y,y_0)\to(\tilde X,\tilde x_0)$ such that
the following diagram commutes:
\begin{center}
\begin{picture}(100,70)(0,-10)
\put(0,0){$(Y,y_0)$}
\put(75,0){$(X,x_0)$}
\put(75,60){$(\tilde X,\tilde x_0)$}

\put(90,55){\vector(0,-1){45}}
\put(92,30){$\rho$}
\put(15,10){\vector(4,3){60}}
\put(43,37){$\tilde \phi$}
\put(30,3){\vector(1,0){42}}
\put(47,-6){$\phi$}
\end{picture}
\end{center}

If $(Y,y_0)$ is ~arcwise connected and locally ~arcwise connected
and $(X,x_0)$ has a covering space $(\tilde X,\tilde x_0,\rho)$,
then by a well-known theorem (c.f. Theorem 5.1 \cite{Ma} p. 128)
a map $\phi :(Y,y_0)\to(X,x_0)$ has a lifting $\tilde\phi:(Y,y_0)
\to (\tilde X,\tilde x_0)$ if and only if $\phi_*(\pi(Y,y_0))
\subseteq\rho_*(\pi(\tilde X,\tilde x_0))$ holds. If such a lifting
exists, it is unique.

\begin{Theorem}
\label{2.1theo}
Let $G$ be a LCA-group with ~arcwise connected dual group $\hat G$ 
and $\C^*=\C-\{0\}$.
Then a probability measure $\mu$ on $G$ is admissible if and
only if the following two conditions hold:
\begin{enumerate}
\item $\hat\mu(\gamma)\not=0\, , ~~\forall\gamma\in\hat G$,
i.e., $\hat\mu$ can be seen as a map $\hat\mu:(\hat G,1)
\to(\C^*,1)$;
\item $\hat\mu_*(\pi(\hat G,1))=\{0\}$ holds for
$\hat\mu:(\hat G,1)\to(\C^*,1)$.
\end{enumerate}

\end{Theorem}
\noindent {\bf Proof.} By Remark \ref{2.1rem} (i), $\hat G$ is
also locally arcwise connected.

Let $R=(0,\infty)\times\R$ be the Riemannian surface of the
logarithm. We shortly denote the element $(1,0)\in R$ by $1$.
We set
$$\rho:(r,\theta)\in R\mapsto re^{i\theta}\in\C ^*$$
then it is clear that $(R,1,\rho)$ is a covering space of
$(\C^*,1)$. Furthermore, we notice that
$$\log_R :(r,\theta)\in R\mapsto\log r+i\theta\in\C$$
maps $(R,1)$ homeomorphically to $(\C,0)$, where the symbol
$\log$ on the right hand side of the mapping stands for the
real logarithm. Thus by using the lifting argument again, we
derive the following commutative diagram:
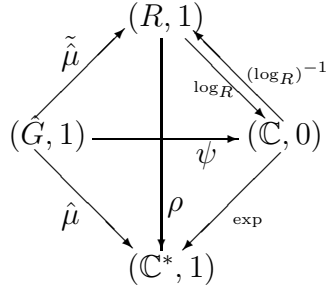
\begin{center}
\begin{picture}(100,100)(0,0)
\put(-0,50){$(\hat G,1)$}
\put(90,50){$(\C,0)$}
\put(45,95){$(R,1)$}
\put(45,0){$(\C^*,1)$}
\put(32,52){\vector(1,0){55}}
\put(71,43){$\psi$}
\put(10,60){\vector(1,1){34}}
\put(20,80){$\tilde {\hat \mu}$}
\put(10,48){\vector(1,-1){37}}
\put(20,20){$\hat \mu$}
\put(58,90){\vector(0,-1){80}}
\put(60,25){$\rho$}
\put(67,90){\vector(1,-1){30}}
\put(90,75){$\scriptscriptstyle (\log_R)^{-1}$}
\put(104,59){\vector(-1,1){34}}
\put(70,70){$\scriptscriptstyle \log_R$}
\put(103,47){\vector(-1,-1){37}}
\put(85,20){$\scriptscriptstyle \exp$}
\end{picture}
\end{center}

Thus, a probability measure $\mu$ on the LCA-group $G$ has a Fourier transform
$\hat \mu$ which can be represented in the form $\hat \mu =e^\psi$ with $\psi:
(\hat G,1)\to(\C,0)$ if and only if $\hat \mu(\gamma)\not=0$,$ 
\forall \gamma\in
\hat G$, and $\hat\mu:(\hat G,1)\to(\C^*,0)$ can be lifted to 
$(R,1)$. Taking into
account the topological properties of $\hat G$ and $\pi(R,1)=0$ (and 
consequently
$\rho_*(\pi(R,1))=\{0\}$), an application the theorem on the 
existence of liftings
given above then concludes the proof.\kasten

\

\begin{Remark}
\label{2.2rem}
{\rm (i) Of course, the second condition of Theorem \ref{2.1theo}
is trivial if $\pi(\hat G,1)=\{0\}$. This is true for $G=\hat G=
\R^d$. In this special case Theorem 1.2 can be obtained without
using homotopy, see e.g. \cite{AGW} for d=1 or \cite{B}
p. 220-223 for $d\in\N$. Therefore, the notion of admissible
probability measures on LCA-groups presented here really extends
the notion of admissibility given in \cite{AGW} in the case that
$G=\R$ (and also extends the discussion presented in \cite{B} for
$G=\R^d, d\in\N$).

\noindent (ii) If the LCA-group $G$ has arcwise connected and
locally arcwise connected dual group $\hat G$ with nontrivial
fundamental group, there exist non admissible probability
measures $\mu\in{\cal M}_1(G)$ with $\hat\mu(\gamma)\not=0
~~\forall\gamma\in\hat G$. Notice that the simplest example for
such a LCA-group $G$ is $\Z$ with dual group $\T^1$, hence let
us take $G=\Z$ as an example to elucidate this point. In this
case, we have
$(\hat G,1)\cong(\T,0)=[0,2\pi]/{\scriptstyle0\sim2\pi}$. Then
$\pi(\hat G,1)\cong\pi(\T,0)\cong\Z$. Furthermore we identify
$\pi(\C^*,1)$ with $\Z$. Let $n\in\Z-\{0\}$ and $\delta_n$ be the
Dirac measure at $n$. For $s\in \T$ we get that $\hat\delta_n(s)
=e^{is n}\not=0$. But $\hat\delta_{n*}:\Z\to\Z$ obviously acts as
multiplication by $n$ and thus $\hat \delta_{n*}(\pi(\hat G,1))
\cong n\Z\not=\{0\}$. Thus, $\delta_n$ is not admissible. Since
a copy of $\T$ is contained in every such $\hat G$, the above
argument carries over to the general case.

\noindent (iii) The above example also shows that for any
LCA-group $G{\subset}V=\R^{d+n}$ with arcwise connected and
locally arcwise connected dual group $\hat G$ (therefore one
has $G{\neq}V$!), in general,
${\cal M}^a_1(G)\not\subseteq{\cal M}_1^a(V)$: We
have for $\mu\in{\cal M}_1(G)\cap{\cal M}^a_1(V)$ with $V$-second
characteristic $\psi$ that $\mu\in{\cal M}^a_1(G)\Leftrightarrow
\psi(x,y)=\psi(x,y+z)\forall x\in\R^d,y\in\R^n, z\in 2\pi\Z^n$.}
\end{Remark}

Remark \ref{2.2rem} (ii) gives rise to the question, ~whether
such non admissible $\mu\in{\cal M}_1(G)$ with $\hat\mu(\gamma)
\not=0 ~ \forall\gamma\in\hat G$ can occur as integer roots of
some admissible probability measure $\nu\in{\cal M}_1^a(G)$?
The following Theorem gives a negative answer.

\begin{Theorem}
\label{2.2theo}
Let $G$ be a LCA-group with locally arcwise and arcwise connected
dual group $\hat G$. For $\mu\in{\cal M}_1^a(G)$ and $n\in\N$, we
assume that $\mu_{1\over n}\in{\cal M}_1(G)$ is a solution to the
problem
\begin{equation}
\label{2.1eqa}
\mu=\underbrace{\mu_{1\over n}\ast\cdots\ast\mu_{1\over n}
}_{n ~ \mbox{\tiny times}}.
\end{equation}
Then $\mu_{1\over n}$ is admissible. Furthermore, $\mu_{1\over n}$
is the only probability measure on $G$ which solves the above
problem.
\end{Theorem}

\noindent{\bf Proof.} The condition that $\hat{\mu_{1\over n}}
(\gamma)\not=0 ~ \forall\gamma\in\hat G$ ~follows immediately by
taking the Fourier transformation of Equation (\ref{2.1eqa}). It
remains to show that condition 2. of Theorem \ref{2.1theo} holds.

Since $\mu_{1\over n}$ solves the above problem, the following
diagram commutes:
\begin{center}
\begin{picture}(100,70)(0,-10)
\put(0,0){$(\hat G,1)$}
\put(75,0){$(\C^*,1)$}
\put(75,60){$(\C^*,1)$}

\put(90,55){\vector(0,-1){45}}
\put(92,30){$z\mapsto z^n$}
\put(15,10){\vector(4,3){60}}
\put(40,43){$\hat{\mu_{1\over n}}$}
\put(30,3){\vector(1,0){42}}
\put(47,-6){$\hat\mu$}
\end{picture}
\end{center}

By identifying $\pi(\C^*,1)$ with $(\Z,0)$, we get that the
following diagram commutes:
\begin{center}
\begin{picture}(100,70)(0,-10)
\put(0,0){$\pi(\hat G,1)$}
\put(75,0){$(\Z,0)$}
\put(75,55){$(\Z,0)$}

\put(90,51){\vector(0,-1){41}}
\put(92,30){$m\mapsto nm$}
\put(15,10){\vector(4,3){58}}
\put(38,43){${\hat{\mu_{{1\over n}}*}}$}
\put(33,3){\vector(1,0){39}}
\put(47,-6){$\hat \mu_*$}
\end{picture}
\end{center}

Since $\hat\mu_*(\pi(\hat G,1))=\{0\}$, this is only possible if
$\hat{\mu_{{1\over n}}*}(\pi(\hat G,1))=\{0\}$. Thus, $\mu_{1\over
n}$ is admissible.

Let $\psi$ be the second characteristic of $\mu$ and $\psi_{1\over
n}$ the second characteristic of $\mu_{1\over n}$. Since $\hat\mu
=(\hat{\mu_{1\over n}})^n=e^{n\psi_{1\over n}}$ and $n\psi_{1\over
n}$ is continuous and fulfils $n\psi_{1\over n}(1)=0$ we get that
$\psi =n\psi_{1\over n}$. But this determines $\psi_{1\over n}$
uniquely and thus $\mu_{1\over n}$ is also uniquely determined.
\kasten

In other words, this theorem says that there exists a solution to
Equation (\ref{2.1eqa}) if and only if ${1\over n}\in\Lambda(\mu)$.
If such a solution exists, then its Fourier transform equals to
$e^{{1\over n}\psi}$, where $\psi$ is the second characteristic of
$\mu$. Since the convolution of an admissible measure with itself
is again admissible, this proves that $\Lambda(\mu)$ really
measures the "divisibility" of $\mu$, i.e. $\Lambda^{\mbox{\rm alg}
}(\mu)=\Lambda(\mu)\cap\Q$.

\begin{Remark}{\rm The uniqueness of the roots can also be obtained
using the imbedding of $G$ into a vector space $V$ and applying the
result of Bauer \cite{B} p.220-223 on the uniqueness of roots of
probability measures on (finite dimensional) vector spaces.
G-admissibility of the roots then follows from Remark \ref{2.2rem}
(iii).}
\end{Remark}
\begin{Remark}
{\rm
The considerations of this section and Theorem \ref{3.1theo} below can be
extended to the case of LCA-groups with only connected (not arcwise connected)
dual group $\hat G$ as follows:

By \cite{He} Theorem N (p. 15), any such $\hat G$ can be obtained as
projective limit of arcwise connected LCA groups $\hat G_{d,n}=\R^d\times\T^n$.

If we assume that the Fourier transform $\hat\mu$ of a probability 
measure $\mu$
on $G$ fulfils the Conditions 1. and 2. of Theorem \ref{2.1theo} if
projected
to any $\hat G_{d,n}$, we get the existence of a second 
characteristic $\psi_{d,n}$
on $\hat G_{d,n}$ and the uniqueness result of Theorem \ref{2.1theo} 
now implies that
the $\psi_{d,n}$ form a projective family with projective limit $\psi$
being the second characteristic defined on $\hat G$. On the other hand,
if a second characteristic $\psi$ exists, then the projection
$\psi_{d,n}$ is defined on $\hat G_{d,n}$ and by Theorem \ref{2.1theo}
the conditions 1. and 2. hold for the projection $\hat \mu_{d,n}$ of
$\hat \mu$ to $\hat G_{d,n}$.
}
\end{Remark}

\section{On the divisibility of non admissible probability
measures}

Let $G$ be a LCA-group with arcwise connected dual group $\hat G$.
  In this section we investigate what happens if for a given probability
measure $\mu$ on $G$ either of the requirements of Theorem
\ref{2.1theo} is not fulfiled. This leads to restrictions on
$\Lambda^{\mbox{alg}}(\mu)$. First we consider this for the case
that condition 1. of Theorem \ref{2.1theo} hold and condition 2.
is violated:

\begin{Theorem}
\label{3.1theo}
Let $\mu\in{\cal M}_1(G)-{\cal M}_1^a(G)$ such that $\hat\mu(
\gamma)\not =0$ holds $\forall\gamma\in\hat G$. We identify
$\pi(\C^*,1)$ with $\Z$ and $n|\hat\mu_*(\pi(\hat G,1))$ means
that $\forall z\in\hat\mu_*(\pi(\hat G,1))\linebreak\exists~ l
\in\Z$ such that $z=nl$. Then
\par\noindent
(i)
$$\Lambda^{\mbox{\rm alg}}(\mu)\subseteq\bigcup_{n\in\N:n|
\hat\mu_*(\pi(\hat G,1))}{1\over n}\N.$$
(ii)
$$\Lambda^{\mbox{\rm alg}}(\mu)\subseteq\bigcap_{z\in\hat\mu_*
(\pi(\hat G,1))-\{0\}}{1\over|z|}\N.$$
(iii) $\Lambda^{\mbox{\rm alg}}(\mu)$ is bounded from below by
$1/\min\{|z|:z\in\hat\mu_*(\pi(\hat G,1))-\{0\}\}$.
\end{Theorem}

\noindent {\bf Proof.} (i) By a diagram analogous to the second
diagram in the proof of Theorem \ref{2.2theo} we get for $q\in
\Lambda^{\mbox{alg}}(\mu), q={m\over n}, m,n\in\N$ mutually
prime that for a solution of Equation (\ref{1.2eqa})
\begin{equation}
\label{3.0eqa}
n~\hat{\mu_{{m\over n}}*}(\pi(\hat G,1))=m~\hat\mu_*(\pi(\hat
G,1))~.
\end{equation}
Since $n$ and $m$ are mutually prime, this can only be true if for
$z\in\hat\mu_*(\pi(\hat G,1))~\exists\linebreak l\in\Z$ s.t.
$z=nl$ holds.

(ii) Let $q,m,n$ be as above. By Eq. (\ref{3.0eqa}) $\forall z\in
\hat \mu_*(\pi(\hat G,1))-\{0\}~ \exists ~ l\in\Z-\{0\}$ such that
$nl=mz$. Thus, $q=m/n=l/z\in{1\over|z|}\N$.

(iii) $nl=mz$ in (ii) implies $n\leq m|z|\Leftrightarrow q\geq1/|z|
~\forall z\in\hat\mu_*(\pi(\hat G,1))-\{0\}$.
\kasten

This shows that the divisibility structure for the class of measures
described in Theorem \ref{3.1theo} {\it a priori} is much simpler
than the one of admissible probability measures. In particular,
there is no open interval $(a,b)\cap\Q$ contained in $\Lambda^{
\mbox{alg}}(\mu)$ for such $\mu$.

Now we consider the situation where the characteristic function
$\hat\mu$ of a non admissible probability measure $\mu\in{\cal
M}_1(G)$ has zero points, i.e. there exists at least one $\gamma_0
\in\hat G$ s.t. $\hat\mu(\gamma_0)=0$.For $\gamma_0\in\hat G$ we
denote by ${\cal P}(\gamma_0)$ the collection of all smooth paths
$\alpha:((-1,1),0)\to (\hat G,\gamma_0)$.

We recall (cf. \cite{BeFo} p. 12) that for any characteristic
function $\hat \mu$ the following inequality holds:
\begin{equation}
\label{3.1eqa}
\left|\hat\mu(\gamma_1)-\hat\mu(\gamma_2)\right|^2\leq2(1-\Re
\hat\mu(\gamma_1\cdot\gamma_2^{-1}))~~\forall\gamma_1,\gamma_2\in
\hat G .
\end{equation}
Let us now assume that $\hat\mu$ is $C^3$-differentiable at $1\in
\hat G$. Note that the real part of $\hat\mu, \Re\hat\mu$, takes
the value $1$ at $1\in\hat G$ and thus takes a maximum at $1\in
\hat G$. For any path $\alpha\in{\cal P}(\gamma_0)$ we then get
by Taylor's formula
$$1-\Re\hat\mu(\alpha(s)\cdot\gamma_0^{-1})=-{1\over 2}\left.{
d^2\Re\hat\mu(\alpha(s)\cdot\gamma_0^{-1})\over ds^2}\right|_{
s=0}s^2+ o(s^2).$$
If we insert this expansion into (\ref{3.1eqa}) and divide by
$s^2$ we get in the limit $s\to+0$
\begin{eqnarray}
\label{3.2eqa}
\limsup_{s\to+0}{|\hat\mu(\alpha(s))|\over s} &=&
\sqrt{\limsup_{s\to +0} {\left| \hat \mu (\alpha(s))\over
s\right|}^2} \nonumber \\
& \leq &\sqrt{-\left. {d^2 \Re \hat\mu (\alpha (s)\cdot
\gamma_0^{-1}) \over ds^2}\right|_{s=0} }<\infty.
\end{eqnarray}

Let ${\cal N}(\hat\mu)=\{\gamma\in\hat G: \hat\mu(\gamma)=0\}
\not =\emptyset$. We define $t_0(\mu)$ by
$$t_0(\mu)=\sup_{\gamma_0\in {\cal N}(\hat\mu)}\sup_{\alpha\in
{\cal P}(\gamma_0)}\sup\{t\geq0:\limsup_{s\to+0}s^{-1}|\hat\mu
(\alpha(s))|^t=\infty\}.$$
Clearly, $t_0(\mu)$ is well-defined, since $\{t\geq0:\limsup_{s\to
+0}s^{-1}|\hat\mu(\alpha(s))|^t=\infty\}$ contains $0$ and is thus
non-empty and is bounded by $1$ (cf. (\ref{3.2eqa})).

We then get
\begin{Theorem}
\label{3.2theo}
Let $\mu$ be a probability measure on $G$ s.t. ${\cal N}(\hat\mu)
\not=\emptyset$. Furthermore, we assume that $\hat\mu$ is
$C^3$-differentiable at $1\in\hat G$.

\noindent Then $\Lambda^{\mbox{\rm alg}}(\mu)$ is bounded from
below by $t_0(\mu)$, i.e. $\Lambda^{\mbox{\rm alg}}(\mu)\subseteq
\Q\cap [t_0(\mu),\infty)$.
\end{Theorem}
{\bf Proof.} Suppose there exists a $t\in\Lambda^{\mbox{\rm alg}}
(\mu)$ with $0<t<t_0(\mu)$. By the definition of $t_0(\mu)$ there
exists a $\gamma_0\in{\cal N}(\hat\mu)$ and a path $\alpha\in{\cal
P}(\gamma_0)$ such that
$$
\limsup_{s\to+0}{|\hat\mu(\alpha(s))|^t\over s}=\infty .
$$
But $t\in\Lambda^{\mbox{\rm alg}}(\mu)$ implies that there exist
$m,n\in \N$ mutually prime with $t={m\over n}$ and a characteristic
function $\hat{\mu_{m\over n}}$ s.t. Equation (\ref{1.2eqa}) holds.
This immediately implies that $|\hat{\mu_{m\over n}}|=|\hat\mu|^t$
holds. Furthermore, the fact that $\hat\mu$ is $C^3$-differentiable
at $1\in\hat G$ implies that $\hat{\mu_{m\over n}}$ is also
$C^3$-differentiable at $1\in \hat G$, since the solution
$\hat{\mu_{n\over m}}$ of the equation (\ref{1.2eqa}) is uniquely
determined in a neighbourhood $U$ of $1\in\hat G$ and depends
smoothly on the values of $\hat\mu$ on $U$. Then by inequality
(\ref{3.2eqa}) one gets that
$$
\limsup_{s\to+0}{|\hat\mu(\alpha(s))|^t\over s}=\limsup_{s\to+0}
{|\hat{\mu_{m\over n}}(\alpha(s))|\over s}<\infty ,$$
which is in contradiction with the assumption that $t\in
\Lambda^{\mbox{\rm alg}}(\mu), ~ 0<t<t_0(\mu)$. \kasten

Theorem \ref{3.2theo} provides an effective tool to calculate
lower bounds of $\Lambda^{\mbox{\rm alg}}(\mu)$ for $\mu\in{\cal
M}_1(G)$ with ${\cal N}(\hat\mu)\not=\emptyset$:
Let e.g. $G=\R$ and $\mu ={1\over2}(\delta_1+\delta_{-1})$. Then
$\hat\mu(y)=\cos y$ has a zero point at $y={\pi\over2}$. Clearly,
$\hat \mu$ is $C^\infty$-differentiable in $0$. Let $\alpha(s)=
{\pi\over 2}+s$ for $s\in(-1,1)$. Since
\begin{eqnarray*}
\limsup_{s\to+0}{|\cos({\pi\over2}+s)|^t\over s} &=& \limsup_{s\to
+0}{|\sin s|^t\over s} \\
&=&\limsup_{s\to+0}|\sin s|^{t-1}=\infty ~ \mbox{for} ~ 0<t<1,
\end{eqnarray*}
we get that $t_0(\mu)\geq1$ and thus $\Lambda^{\mbox{\rm alg}}(
\mu)\subseteq\Q\cap[1,\infty)$. Since $1\in\Lambda^{\mbox{\rm alg}
}(\mu)$ this estimate is sharp.

\section{The case of LCA-groups with non connected dual groups}
In this section we give a short discussion, why one can not
expect results similar to those of Section 2 and Theorem
\ref{3.1theo} for the case of a LCA-group
with non-arcwise connected dual group. We do not intend treat
partly divisibility for this case exhaustively, but we want to
point out that in this case the divisibility of a probability
measure $\mu$ with a $\hat \mu =e^\psi$ in general is not
associated to a closed semigroup $\Lambda(\mu)$.

As the simplest example for a LCA-group with non connected dual
group we study the case $G=\Z_2$. Any probability measure on $G$
is of the form $\mu =\alpha\delta_0+(1-\alpha)\delta_1$, $\alpha
\in [0,1]$. Identifying $\hat G$ with $\Z_2$ we get for the
Fourier transform of $\mu$ that $\hat\mu (0)=1, ~ \hat\mu(1)=2
\alpha-1$. If $\alpha\in(1/2,1]$, then $\psi=\log\hat\mu$ is well
defined and $e^{t\psi}$ for $t>0$ is the Fourier transform of the
probability measure $\mu_t=\alpha_t\delta_0+(1-\alpha_t)\delta_1,~
\alpha_t=((2\alpha-1)^t+1)/2$. Thus, $\mu$ is infinitely divisible.
This consideration extends to the case $\alpha=1/2$ if we adopt
the convention $\log0=-\infty, e^{-\infty}=0$.

For $\alpha\in(0,1/2)$ however, we have the following well known
result:

\begin{Theorem}
\label{4.1theo}
Let $\alpha\in(0,1/2)$ and $\mu=\alpha\delta_0+(1-\alpha)\delta_1$
be given. Then the following equation
$$\mu=\underbrace{\mu_{1\over n}\ast\cdots\ast\mu_{1\over n}}_{n
~ \mbox{\tiny times}}  \,  , n\in\N
$$
has a solution $\mu_{1\over n}$ if and only if $n$ is odd.
\end{Theorem}

\noindent{\bf Proof.} The if part follows from $\mu=\delta_1\ast
\nu$ with $\nu=(1-\alpha)\delta_0+\alpha\delta_1$ since $\nu$ is
infinitely divisible and $\delta_1$ convoluted $n$ times with
itself gives $\delta_1$ if $n$ is odd.

For the only if part suppose $\mu_{1\over n}$ solves the above
equation, $\mu_{1\over n}=\beta\delta_0+(1-\beta)\delta_1, \beta
\in [0,1]$. But then $(\hat{\mu_{1\over n}}(1))^n=(2\beta-1)^n=
2\alpha-1<0$ leads to a contradiction if $n$ is even. \kasten

Theorem \ref{4.1theo} is interesting in this context, since it
makes clear that for $\alpha\in(0,1/2)$ the set $\Lambda^{\mbox{
alg}}(\mu)$ can not be imbedded into a closed semigroup $S\subseteq
(0,\infty)$ such that $\Lambda^{\mbox{alg}}(\mu)=S\cap\Q$. This
follows e.g. from Remark \ref{2.1rem}
(iii).

Similar effects occur in the case when $\Z_2$ is replaced by $\Z_n$ (or
any other discrete group $G$ with nilpotent elements):
If we consider e.g. the measure $\delta_1$ on $\Z_n$ it is immediately
clear that $\Lambda^{\mbox{alg}}(\delta_1)=
\{m/l:m,l\in\N\mbox{ m.p. } m=l\mbox{ mod }n\}\not=(0,\infty)$ and contains
sequences converging to zero, and thus by Remark \ref{2.1rem} (iii), again
cannot be imbedded into a closed semigroup.

The reason that the situation described here is different from that
in the previous sections is that the second characteristic $\psi$
associated to a probability measure $\mu$ is no longer unique in the
case that $\hat G$ is not connected. More explicitly, for any second
characteristic $\psi$, we have that $\psi_k=\psi+2\pi ik$ with
$k:\hat G\to\Z$ is constant on the connected components of $\hat G$
and fulfilling $k(1)=0$ is another continuous function with $\hat
\mu=e^{\psi_k}$ and $\psi_k(1)=0$.

Thus, $\Lambda^{\mbox{alg}}(\mu)$ is the rational part of an
(infinite) union of closed semigroups
$$\Lambda_k(\mu)=\{t>0: e^{t\psi_k}\mbox{is a characteristic
function}\}.$$
But such a union in general is neither closed nor a semigroup. This
shows that for non connected $\hat G$ the set $\Lambda^{\mbox{alg}}
(\mu)$ loses a lot of its structure. For this reason it is much less
attractive to study this case than the case where $\hat G$ is
connected.

\begin{Remark}
{\rm
(i) The non-uniqueness of logarithms for infinitely divisible measures on
(not necessarily abelian) finite groups has been studied systematically by
J. B\"oge \cite{Bo}. Even though the considerations of \cite{Bo} do not have
a straight forward extension to the case of partly divisible measures (since
it uses a compactness argument for the construction of logarithms which can
not be applied if the roots are in the non-compact set of complex measures),
the non-uniqueness of the representation $\delta_0=\exp\psi$ where $\exp$ is
defined on the group algebra and $\exp(t\psi)$ is a measure for $t>0$, in
general creates a non-uniqueness of roots also in the partial divisible case:
If $\mu_{1\over n}$ is an $n$-th root of $\mu$, then this is also true for
$\mu_{1\over n}'=\mu{1\over n}*\exp(\psi/n)$ for some $\psi\not=0$ and we
expect that in general $\mu_{1\over n}\not=\mu_{1\over n}'$.

\

\noindent (ii) The nilpotency of elements also leads to the non-uniqueness
of roots on LCA-groups containing a torus. E.g. on $\T$ there are exactly $n$
$n$-th roots of $\delta_0$. The situation is a little bit different from the
case of finite groups, since these roots are infinitely divisible and thus do
not give examples for measures $\mu$ with $\Lambda^{\mbox{alg}}(\mu)$ not
embeddable in a closed semigroup $S\subseteq (0,\infty)$. But still there is
no generalisation of Theorem \ref{2.2theo} also for this case.
}
\end{Remark}

\bigskip\noindent
{\bf Acknowledgements.}
Stimulating discussions with Dr. W. R. Schneider and the financial
support of D.F.G. (SFB 237) are gratefully acknowledged. We also
thank the referees for their careful reading of the paper and their
useful remarks and suggestions.

\bigskip\medskip\par\noindent
Institut f\"ur Angewandte Mathematik der Universit\"at Bonn,
\par\noindent
\quad D-53013 Bonn, Germany
\par\noindent
SFB 256 Bonn, Germany
\par\noindent
SFB 237 Essen-Bochum-D\"usseldorf, Germany
\par\noindent
BiBoS Research Centre, Bielefeld-Bochum, Germany
\par\noindent
CERFIM, Locarno, Switzerland
\par\noindent
Acc. Arch., Mendrisio, Switzerland

\bigskip\medskip\par\noindent
Dipartimento di Matematica, Universit\`a di Roma "La Sapienza",
\par\noindent
\quad Piazzale Aldo Moro 2, I-00185 Roma, Italy

\bigskip\medskip\par\noindent
Fakult\"at f\"ur Mathematik der Ruhr-Universit\"at,
D-44780 Bochum, Germany
\par\noindent
SFB 237 Essen-Bochum-D\"usseldorf, Germany
\par\noindent
Institute of Applied Mathematics, Academia Sinica,
Beijing 100080, P R China


\begin{thebibliography}{99}

\bibitem{AGW} Albeverio, S., Gottschalk, H., Wu, J.-L.: Partly
divisible probability distributions. Forum Math. {\bf 10}, 687--697
(1998).

\bibitem{B} Bauer, H.:Wahrscheinlichkeitstheorie, 4th ed.
Berlin: de Gruyter 1991.

\bibitem{BeFo} Berg, C., Forst, G.: Potential theory on locally
compact Abelian groups. Berlin: Springer-Verlag
1975.

\bibitem{BS} Bisgaard, T.M., Sasv\'ari, Z.: On the positive
definiteness of certain functions. Math. Nachr. {\bf 186}, 81-99
(1997).

\bibitem{Bo} B\"oge, J.: \"Uber die Charakterisierung unendlich
teilbarer Wahrscheinlichkeitsverteilungen. J. Reine Angew. Math.
{\bf 201}, 150-156 (1959).

\bibitem{Di} Dixmier, M. J.:Quelques propri\'et\'es des groupes
ab\'eliens localement compacts. Bull Sci. Math {\bf 81}, 38-48
(1957).

\bibitem{Du} Dugu\'e, D.:
Arithm\'etique des lois de probabilit\'es. M\'emoriale des Sciences
Math., Vol. {\bf 137}, Paris: Gauthier-Villars 1957.

\bibitem{GW} Gillespie, T.A., West, T. T.: Weakly compact groups of
operators. Proc. AMS {\bf 49}, No. 1, 78-83 (1975).

\bibitem{Ha} Hazod, W.: \"Uber Wurzeln und Logarithmen
beschr\"ankter Ma\ss e. Z. Wahrscheinlichkeitstheorie verw. Geb.
{\bf 20}, 259-270 (1971).

\bibitem{He} Heyer, H.: Probability Measures on Locally Compact
Groups. Berlin: Springer-Verlag 1977.

\bibitem{Hi} Hille, E.: Roots and logarithms of elements of a
complex Banach algebra. Math. Annalen {\bf 136}, 46-57 (1958).

\bibitem{Ma} Massay, W.S.: A basic course in algebraic topology.
New York: Springer-Verlag 1991.

\end{thebibliography}
\end{document}